\documentclass[10pt,leqno,a4paper]{amsart}

\newtheorem{lemma}{Lemma}
\newtheorem{theorem}{Theorem}
\newtheorem{corollary}{Corollary}
\newtheorem{proposition}{Proposition}
\newtheorem{definition}{Definition}

\newtheorem{example}{Example}
\newtheorem{remark}{Remark}
\usepackage{graphics}
\usepackage{epsfig}
\usepackage{amssymb,amsfonts,amscd}
\usepackage{enumerate}
\def\C{{\mathbb C}}

\def\R{{\mathbb R}}
\def\B{{\mathbb B}}
\def\CC{{\mathcal C}}
\begin{document}

\title[Estimates of the Kobayashi-Royden metric]
{Estimates of the Kobayashi-Royden metric in almost complex manifolds}

\author{Herv\'e Gaussier and Alexandre Sukhov}

\address{\begin{tabular}{lll}
Herv\'e Gaussier & & Alexandre Sukhov\\
C.M.I. & &U.S.T.L. \\
39, rue Joliot-Curie,& &Cit\'e Scientifique \\
13453 Marseille Cedex 13 & &59655 Villeneuve d'Ascq Cedex
\end{tabular}
}

\email{gaussier@cmi.univ-mrs.fr, \ sukhov@agat.univ-lille1.fr} 


\date{\number\year-\number\month-\number\day}

\begin{abstract}
We establish a lower estimate for the Kobayashi-Royden infinitesimal
pseudometric on an almost complex manifold $(M,J)$ admitting a bounded
strictly plurisubharmonic function. We apply this result to study the
boundary behaviour of the metric on a strictly pseudoconvex domain in
$M$ and to give a sufficient condition for the complete hyperbolicity
of a domain in $(M,J)$. 
\end{abstract}

\maketitle
\section*{Introduction}
In the recent paper \cite{ko01}, S.Kobayashi studied the following question~:
Does every point in an almost complex manifold admit a basis of
complete hyperbolic neighborhoods ? 
This question was solved in dimension four by R.Debalme and
S.Ivashkovich in \cite{di01}.

In the present paper, we give a lower estimate
on the Kobayashi-Royden infinitesimal metric on a strictly pseudoconvex
domain in an almost complex manifold (such estimates are well-known in
the integrable case  \cite{gr75}). A corollary of our main resultgives a
positive answer to the previous question in any dimension~:
~every point in an almost complex manifold has a complete hyperbolic
neighborhood.

Our approach consists of two parts.  In order to localize the
Kobayashi-Royden metric we use a method developped essentially by
N.~Sibony \cite{si81} in the case of the standard complex structure
(see also \cite{df77} by K.~Diederich-J.E.~Fornaess and \cite{kr81} by
N.~Kerzman-J.P.~Rosay). This is based on the construction of special
classes of plurisubharmonic functions.  Then we apply an almost
complex analogue of the scaling method due to S.Pinchuk in the
integrable case (see, for instance, \cite{pi91}) and obtain presice
estimates of the metric.  We point out that similar ideas have been
used by F.Berteloot \cite{be95,be02} in order to estimate the
Kobayashi-Royden metric on some classes of domains in $\mathbb C^n$.
\vskip 0,1cm
We note that S.Ivashkovich and J.P.Rosay recently proved
in \cite{iv-ro04}, among other results, estimates of the Kobayashi-Royden
metric similar to ours under weaker assumptions on the regularity of the
almost complex structure.

\vskip 0,1cm
\noindent {\sl Acknowledgments}. The authors thank E.~Chirka,
 B.~Coupet, S.~Ivashkovich and J.-P.Rosay for helpful discussions and
 the referee for valuable remarks. We are particularly indebted to
S.Ivashkovich who pointed out an erroneous argument in the previous version
of our paper.

\section{Preliminaries}

\subsection{Almost complex manifolds.} 
 Let $(M',J')$ and $(M,J)$ be almost complex manifolds and let $f$ be
 a smooth map from $M'$ to $M$. We say that $f$ is {\sl
 $(J',J)$-holomorphic} if $df \circ J' = J \circ df$ on $TM'$. 
We denote by $\mathcal O_{(J',J)}(M',M)$ the set of
$(J',J)$-holomorphic maps from $M'$ to $M$. Let $\Delta$ be the unit
disc in $\C$ and $J_{st}$ be the standard integrable structure on $\C^n$
for every $n$. If $(M',J')=(\Delta,J_{st})$,
we denote by $\mathcal O_J(\Delta,M)$ the set
$\mathcal O_{(J_{st},J)}(\Delta,M)$ of {\sl $J$-holomorphic discs} in $M$. 

The following Lemma shows that every almost complex manifold
$(M,J)$ can be viewed locally as the unit ball in
$\mathbb C^n$ equipped with a small almost complex
deformation of $J_{st}$. This will be used frequently in the sequel.
\begin{lemma}
\label{suplem1}
Let $(M,J)$ be an almost complex manifold. Then for every point $p \in
M$ and every $\lambda_0 > 0$ there exist a neighborhood $U$ of $p$ and a
coordinate diffeomorphism $z: U \rightarrow \mathbb B$ such that
$z(p) = 0$, $dz(p) \circ J(p) \circ dz^{-1}(0) = J_{st}$  and the
direct image $z_*(J) := dz \circ J \circ dz^{-1}$ satisfies
$\vert\vert z_*(J) - J_{st}
\vert\vert_{\CC^2(\bar {\mathbb B})} \leq \lambda_0$.
\end{lemma}
\proof There exists a diffeomorphism $z$ from a neighborhood $U'$ of
$p \in M$ onto $\mathbb B$ satisfying $z(p) = 0$ and $dz(p) \circ J(p)
\circ dz^{-1}(0) = J_{st}$. For $\lambda > 0$ consider the dilation
$d_{\lambda}: t \mapsto \lambda^{-1}t$ in $\C^n$ and the composition
$z_{\lambda} = d_{\lambda} \circ z$. Then $\lim_{\lambda \rightarrow
0} \vert\vert (z_{\lambda})_{*}(J) - J_{st} \vert\vert_{\CC^2(\bar
{\mathbb B})} = 0$. Setting $U = z^{-1}_{\lambda}(\mathbb B)$ for
$\lambda > 0$ small enough, we obtain the desired statement. \qed

\subsection{ $\partial_J$ and $\bar{\partial}_J$ operators}

Let $(M,J)$ be an almost complex manifold. We denote by $TM$ the real 
tangent bundle of $M$ and by $T_\C M$ its complexification. Recall
that $T_\C M = T^{(1,0)}M \oplus T^{(0,1)}M$ where
$T^{(1,0)}M:=\{ X \in T_\C M : JX=iX\} = \{\zeta -iJ \zeta, \zeta \in
TM\},$ 
and $T^{(0,1)}M:=\{ X \in T_\C M : JX=-iX\} = \{\zeta +
iJ \zeta, \zeta \in TM\}$.
 Let $T^*M$ denote the cotangent bundle of  $M$.
Identifying $\C \otimes T^*M$ with
$T_\C^*M:=Hom(T_\C M,\C)$ we define the set of complex
forms of type $(1,0)$ on $M$ by~:
$
T_{(1,0)}M=\{w \in T_\C^* M : w(X) = 0, \forall X \in T^{(0,1)}M\}
$
and the set of complex forms of type $(0,1)$ on $M$ by~:
$
T_{(0,1)}M=\{w \in T_\C^* M : w(X) = 0, \forall X \in T^{(1,0)}M\}.
$.
Then $T_\C^*M=T_{(1,0)}M \oplus T_{(0,1)}M$.   
This allows to define the operators $\partial_J$ and
$\bar{\partial}_J$ on the space of smooth functions defined on
$M$~: given a complex smooth function $u$ on $M$, we set $\partial_J u =
du_{(1,0)} \in T_{(1,0)}M$ and $\bar{\partial}_Ju = du_{(0,1)}
\in T_{(0,1)}M$. As usual,
differential forms of any bidegree $(p,q)$ on $(M,J)$ are defined
by means of the exterior product.

\section{$J$-plurisubharmonic functions with 
logarithmic singularities}

\subsection{Plurisubharmonic functions.} We first recall the following 
definition~:
\begin{definition}\label{d6}
An upper semicontinuous function $u$ on $(M,J)$ is called 
{\sl $J$-plurisubharmonic} on $M$ if the composition $u \circ f$ 
is subharmonic on $\Delta$ for every $f \in \mathcal O_J(\Delta,M)$.
\end{definition} 
If $M$ is a domain in $\C^n$ and $J=J_{st}$ then a 
$J_{st}$-plurisubharmonic function is a plurisubharmonic function 
in the usual sense. 

\begin{definition} Let $u$ be a $\CC^2$ function on $M$, let $p \in M$
and $v \in T_pM$. The Levi 
form of $u$ at $p$, evaluated on $v$, is defined by
$\mathcal L^J(u)(p)(v):=-d(J^\star du)(X,JX)(p)$
where $X$ is any vector field on $TM$ such that $X(p) = v$.
\end{definition}
Following \cite{de99,ha02} we have~:
\begin{proposition}\label{PROP}
Let $u$ be a $\CC^2$ real valued function on $M$, let $p \in M$ and $v \in
T_pM$. Then $\mathcal L^J(u)(p)(v) = \Delta(u \circ f)(0)$
where $f$ is any $J$-holomorphic disc in $M$ satisfying $f(0) = p,\ 
df(0)(\partial / \partial x) = v$.
\end{proposition}

Obviously the Levi form is invariant with respect to biholomorphisms.
More precisely let $u$ be a $\CC^2$ real valued function on $M$,
let $p \in M$ and $v \in T_pM$.
If $\Phi$ is a diffeomorphism from $(M,J)$ to $(M',J')$,
$(J,J')$-holomorphic, then
$\mathcal L^J(u)(p)(v) =
\mathcal L^{J'}(u \circ \Phi^{-1})(\Phi(p))(d\Phi(p)(v))$.

Finally it follows from Proposition~\ref{PROP} that a $\CC^2$ real valued
function $u$ on $M$ is
$J$-plurisubharmonic on $M$ if and only if $\mathcal L^J(u)(p)(v) \geq 0$
for every $p \in M$, $v \in T_pM$.
\vskip 0,1cm
This leads to the definition~:
\begin{definition}
A $\CC^2$ real valued function $u$ on $M$ is {\sl strictly 
$J$-plurisubharmonic} on $M$ if  $\mathcal L^J(u)(p)(v)$
is positive for every $p \in M$, $v \in T_pM \backslash \{0\}$.
\end{definition}

We have the following example of a
$J$-plurisubharmonic function on an almost complex manifold $(M,J)$~:
\begin{example}\label{example}
For every point $p\in (M,J)$ there exists a neighborhood $U$ of $p$
 and a diffeomorphism $z:U \rightarrow \mathbb
B$ centered at $p$ (ie $z(p) =0$) such that the function $|z|^2$ is
$J$-plurisubharmonic on $U$. 
\end{example}
\proof Let $p \in M$, $U_0$ be a neighborhood of $p$ and $z: U_0
\rightarrow \mathbb B$ be local complex coordinates centered at $p$,
such that $dz \circ J(p) \circ dz^{-1} = J_{st}$ on $\mathbb B$. Consider
the function $u(q) = |z(q)|^2$ on $U_0$. For every $w,v \in \C^n$ we
have $\mathcal L^{J_{st}}(u\circ z^{-1})(w)(v) = \|v\|^2$. Let
$B(0,1/2)$ be the ball centered at the origin with radius $1/2$ and
let $\mathcal E$ be the space of smooth almost complex structures
defined in a neighborhood of $\overline{B(0,1/2)}$. Since the function
$(J',w) \mapsto \mathcal L^{J'}(u \circ z^{-1})(w)$ is continuous on
$\mathcal E \times B(0,1/2)$, there exist a neighborhood $V$ of the
origin and positive constants $\lambda_0$ and $c$ such that $\mathcal
L^{J'}(u\circ z^{-1})(q)(v)\geq c \|v\|^2$ for every $q \in V$ and
for every almost complex structure $J'$ satisfying
$\|J'-J_{st}\|_{\CC^2(\bar{V})} \leq \lambda_0$. Let $U_1$ be a
neighborhood of $p$ such that $\|z_*(J) -
J_{st}\|_{\CC^2(\overline{z(U_1)})}\leq \lambda_0$ and let $0<r<1$ be
such that $B(0,r) \subset V$ and $U:=z^{-1}(B(0,r)) \subset U_1$. Then
we have the following estimate for every $q \in U$ and $v \in T_qM$~:
$\mathcal L^J(u)(q)(v) \geq c \|v\|^2$. Then $r^{-1}z$ is the desired
diffeomorphism. \qed

\vskip 0,1cm     
We also have the following 

\begin{lemma}
A function  $u$ of class $\CC^2$ in a neighborhood
of a point $p$ of $(M,J)$  is strictly $J$-plurisubharmonic
if and only there exists a neighborhood $U$ of $p$  with local
complex coordinates $z:U \rightarrow \mathbb B$ centered at $p$, such
that the function $u - c|z|^2$ is $J$-plurisubharmonic on $U$ for some
constant $c > 0$.
\end{lemma}

The function $\log \vert z \vert$ is
$J_{st}$-plurisubharmonic on $\C^n$ and plays an important role in the
pluripotential theory as the Green function for the complex
Monge-Amp\`ere operator on the unit ball. In particular, this function
is crucially used in Sibony's method  in order to localize and
estimate the Kobayashi-Royden metric on a complex manifold. Unfortunately,
after an arbirarily small general almost complex deformation of the
standard structure this function is {\it not} plurisubharmonic with
respect to the new structure (in any neighborhood of the origin), see
for instance \cite{de99}. So we will need the following statement
communicated to the authors by E.Chirka~:
\begin{lemma}
Let $p$ be a point in an almost complex manifold $(M,J)$. There exist
a neighborhood $U$ of $p$ in $M$, a diffeomorphism $z : U \rightarrow
\B$ centered at $p$ and positive constants $\lambda_0,\ A$, such that
the function $\log|z| + A|z|$ is
$J'$-plurisubharmonic on $U$ for every almost complex structure $J'$
satisfying $\|J'-J\|_{\CC^2(\bar{U})} \leq \lambda_0$.
\end{lemma}

\noindent{\it Proof.} Consider the function $u=|z|$ on $\mathbb B$.
Since $\mathcal L^{J_{st}}(u \circ z^{-1})(w)(v) \geq \|v\|^2/4|w|$
for every $w \in \B\backslash \{0\}$ and every $v \in \C^n$, it
follows by a direct expansion of ${\mathcal L}^{J'}(u)$ that there
exist a neighborhood $U$ of $p$, $U \subset\subset U_0$, and a
positive constant $\lambda_0$ such that ${\mathcal L}^{J'}(u)(q)(v)
\geq \|v\|^2/5|z(q)|$ for every $q \in U \backslash\{p\}$, every $v\in
T_qM$ and every almost complex structure $J'$ satisfying
$\|J'-J\|_{\CC^2(\bar{U})} \leq \lambda_0$. Moreover, computing the
Laplacian of $\log|f|$ where $f$ is any $J$-holomorphic disc we
obtain, decreasing $\lambda_0$ if necessary, that there exists a
positive constant $B$ such that $\mathcal L^{J'}(\log|z|)(q)(v) \geq
-B\|v\|^2/|z(q)|$ for every $q \in U\backslash\{p\}$, every $v \in
T_qM$ and every almost complex structure $J'$ satisfying
$\|J'-J\|_{\CC^2(\bar{U})} \leq \lambda_0$. We may choose $A = 2B$ to
get the result. \qed

\section{Localization of the Kobayashi-Royden metric on almost complex
manifolds}

Let $(M,J)$ be an almost complex
manifold. In what follows we use the notation $\zeta=x+iy \in \mathbb C$.
According to \cite{nw63}, for every $p \in M$ there is a
neighborhood $\mathcal V$ of $0$ in $T_pM$ such that for every $v \in
\mathcal V$ there exists $f \in \mathcal O_J(\Delta,M)$ satisfying
$f(0) = p,$ $df(0) (\partial / \partial x) = v$. This allows to define
the Kobayashi-Royden infinitesimal pseudometric $K_{(M,J)}$.
\begin{definition}\label{dd}
For $p \in M$ and $v \in T_pM$, $K_{(M,J)}(p,v)$ is the infimum of the
set of positive $\alpha$ such that there exists a $J$-holomorphic disc
$f:\Delta \rightarrow M$ satisfying $f(0) = p$ and $df(0)(\partial
/\partial x) = v/\alpha$.
\end{definition}
Since for every $f \in \mathcal O_{(J',J)}(M',M)$ and
every $\varphi \in \mathcal O_J(\Delta,M')$ the composition $f \circ
\varphi$ is in $\mathcal O_J(\Delta,M)$ we have~:
\begin{proposition}\label{ppp}
Let $f:(M',J') \rightarrow (M,J)$ be a $(J',J)$-holomorphic map. 
Then $K_{(M,J)}(f(p'),df(p')(v')) \leq K_{(M',J')}(p',v')$for every
$p'\in M', \ v' \in T_{p'}M'$.
\end{proposition}
We denote by $d_{(M,J)}^K$ the integrated pseudodistance of the
Kobayashi-Royden infinitesimal pseudometric. According to the almost
complex version of Royden's theorem \cite{kr99}, it coincides with the
usual Kobayashi pseudodistance on $(M,J)$ defined by means of
$J$-holomorphic discs.  Similarly to the case of the integrable
structure we have~:
\begin{definition}\label{dddd} 
$(i)$ Let $p \in M$. Then $M$ is {\rm locally hyperbolic at $p$} if 
there exists a neighborhood $U$ of $p$ and a positive constant $C$ such 
that for every $q \in U$, $v \in T_qM$~: $K_{(M,J)}(q,v) \geq C \|v\|$.

$(ii)$ $(M,J)$ is {\rm hyperbolic} if it is locally hyperbolic
at every point.

$(iii)$ $(M,J)$ is {\rm complete hyperbolic} if the Kobayashi 
ball $B_{(M,J)}^K(p,r):=\{q \in M : d_{(M,J)}^K(p,q) < r\}$ is relatively 
compact in $M$ for every $p \in M$, $r \geq 0$.
\end{definition}

\begin{lemma}\label{lemlem}
Let $r < 1$ and let $\theta_r$ be a 
smooth nondecreasing function on
$\R^+$ such that $\theta_r(s)= s$ for $s \leq r/3$ and $\theta_r(s) =
1$ for $s \geq 2r/3$. Let $(M,J)$ be an almost complex manifold, and
let $p$ be a point of $M$. Then there exists a neighborhood $U$ of
$p$, positive constants $A = A(r)$, $B=B(r)$ and a diffeomorphism $z:U
\rightarrow \mathbb B$ such that $z(p) = 0$, $dz(p) \circ J(p) \circ
dz^{-1}(0) = J_0$ and the function ${\rm log}(\theta_r(\vert z
\vert^2)) + \theta_r(A\vert z \vert) + B\vert z \vert^2$ is
$J$-plurisubharmonic on $U$.
\end{lemma}

\noindent{\sl Proof of Lemma~\ref{lemlem}}. 
Denote by $w$ the standard coordinates in $\C^n$. It follows
from Lemma~3 that there exist positive constants $A$ and $\lambda_0$
such that the function ${\rm log}(\vert w \vert^2) + A \vert w \vert$
is $J'$-plurisubharmonic on $\mathbb B$ for every almost complex
structure $J'$, defined in a neighborhood of $\bar{\B}$ in $\C^n$ and
such that $\|J'-J_{st}\|_{\CC^2(\bar{\mathbb B})} \leq \lambda_0$. This
means that the function $v(w) = \log(\theta_r(|w|^2)) + \theta_r(A|w|)$
is $J'$-plurisubharmonic on $B(0,r')=\{w \in \C^n : |w| < r'\}$ for
every such almost complex structure $J'$, where $r'=inf(\sqrt{r/3},
r/3A)$. Decreasing $\lambda_0$ if necessary, we may assume that the
function $|w|^2$ is strictly $J'$-plurisubharmonic on $\B$. Then,
since $v$ is smooth on $\mathbb B \backslash B(0,r')$, there exists a
positive constant $B$ such that the function $v + B\vert w \vert^2$ is
$J'$-plurisubharmonic on $\mathbb B$ for
$\|J'-J_{st}\|_{\CC^2(\bar{\mathbb B})} \leq \lambda_0$. According to
Lemma \ref{suplem1} there exists a neighborhood $U$ of $p$ and a
diffeomorphism $z:U \rightarrow \mathbb B$ such that $\vert\vert
z_*(J) - J_{st} \vert\vert_{\mathcal C^2(\bar{\mathbb B})} \leq
\lambda_0$. Then the function $v \circ z = {\rm
log}(\theta_r(\vert z \vert^2)) + \theta_r(A \vert z \vert) +
B\vert z \vert^2$
is $J$-plurisubharmonic on $U$. \qed

\begin{proposition}\label{thm3}
(Localization principle) Let $D$ be a domain in an almost complex
manifold $(M,J)$, let $p \in \bar{D}$, let $U$ be a neighborhood of
$p$ in $M$ (not necessarily contained in $D$) and let $z:U \rightarrow
\mathbb B$ be the diffeomorphism given by Lemma~\ref{lemlem}.  
Let $u$ be a $\mathcal C^2$ function on $\bar{D}$, negative and
$J$-plurisubharmonic on $D$. We assume that $-L \leq u < 0$ on $D \cap
U$ and that $u-c|z|^2$ is $J$-plurisubharmonic on $D \cap U$, where
$c$ and $L$ are positive constants. Then there exist a positive
constant $ s$ and a neighborhood $V \subset \subset U$ of $p$,
depending on $c$ and $L$ only, such that for $q \in D \cap V$ and $v
\in T_qM$ we have the following inequality~:

\begin{equation}\label{e2}
K_{(D,J)}(q,v) \geq s K_{(D \cap U,J)}(q,v).
\end{equation}
\end{proposition}

We note that a similar statement was 
obtained by F.Berteloot \cite{be95} in the integrable case. The proof
is based on N.Sibony's method \cite{si81}.

\vskip 0,2cm
\noindent{\sl Proof of Proposition~\ref{thm3}}.
 Let $0<r<1$ be such that the set $V_1:=\{q \in U :
|z(q)| \leq \sqrt{r}\}$ is relatively compact in $U$ and let $\theta_r$ 
be a smooth nondecreasing function on $\R^+$ such that $\theta_r(s)= s$ for
$s \leq r/3$ and $\theta_r(s) = 1$ for $s \geq 2r/3$. According to Lemma
~\ref{lemlem}, there exist uniform positive constants $A$ and $B$ such
that the function 
$$
{\rm log}(\theta_r(|z-z(q)|^2))+ \theta_r(A|z-z(q)|)+ B|z|^2
$$
is $J$-plurisubharmonic on $U$ for every $q \in V$. 
By assumption the function
$u-c|z|^2$ is $J$-plurisubharmonic on $D \cap U$. Set $\tau=2B/c$ and
define, for every point $q \in V$, the function $\Psi_{q}$ by~:
$$
\left\{
\begin{array}{lll}
\Psi_{q}(z) &=& \theta_r(|z-z(q)|^2)\exp(\theta_r(A|z-z(q)|)) 
\exp(\tau u(z))\ {\rm if} \  z \in D \cap U,\\
& & \\
\Psi_{q} &=& \exp(1+\tau u) \ {\rm on} \ D \backslash U.
\end{array}
\right.
$$

Then for every $0 < \varepsilon \leq B$, the function ${\rm
log}(\Psi_{q})-\varepsilon|z|^2$ is $J$-plurisubharmonic on $D \cap U$
and hence $\Psi_{q}$ is $J$-plurisubharmonic on $D \cap U$. Since
$\Psi_{q}$ coincides with $\exp(\tau u)$ outside $U$, it is globally
$J$-plurisubharmonic on $D$. 

Let $f \in \mathcal O_{J}(\Delta,D)$ be such that $f(0)=q \in V_1$ and
$(\partial f/\partial x)(0) = v/\alpha$ where $v \in T_qM$ and $\alpha
>0$. For $\zeta$ sufficiently close to 0 we have
$$
f(\zeta) = q + df(0)(\zeta) +
\mathcal O(|\zeta|^2).
$$
Setting $\zeta= \zeta_1+i\zeta_2$ and using
the $J$-holomorphy condition $df(0)\circ J_{st} = J \circ
df(0)$, we may write~:
$$
df(0)(\zeta) = \zeta_1 df(0)(\partial /
\partial x) + \zeta_2 J(df(0)(\partial / \partial x)).
$$
Consider the function
$$
\varphi(\zeta) = \Psi_q(f(\zeta))/|\zeta|^2
$$
which is subharmonic on
$\Delta \backslash \{0\}$. Since
$$
\varphi(\zeta) = |f(\zeta)-q|^2/|\zeta|^2 \exp(A|f(\zeta)-q|) 
\exp(\tau u(f(\zeta)))
$$
for $\zeta$ close to 0 and 
$$
|df(0)(\zeta)| \leq |\zeta| (\|I+J\|\,\|df(0)(\partial
/\partial x)\|)
$$
we obtain that $\limsup_{\zeta \rightarrow 0}\varphi(\zeta)
$ is finite. Moreover setting $\zeta_2=0$ we have 
$$
\limsup_{\zeta \rightarrow 0}\varphi(\zeta) \geq \|df(0)(\partial
/\partial x)\|^2\exp(-2B|u(q)|/c).
$$
Applying the maximum principle to a subharmonic extension of $\varphi$
on $\Delta$ we obtain the inequality 
$$
\|df(0)(\partial / \partial x)\|^2 \leq \exp(1+2B|u(q)|/c).
$$
Hence, by definition of the Kobayashi-Royden infinitesimal pseudometric, 
we obtain for every $q \in D \cap V_1$, $v \in T_qM$~:
\begin{eqnarray}
\label{localhyp}
K_{(D,J)}(q,v) \geq \left(\exp\left(-1-2B\frac{|u(q)|}{c}\right)
\right)^{1/2}\|v\|.
\end{eqnarray}
Consider now the Kobayashi ball $B_{(D,J)}(q,\alpha)=\{w \in D :
d_{(D,J)}^K(w,q)<\alpha\}$. It follows from Lemma~2.2 of \cite{ccs99} (whose
proof is identical in the almost complex setting) that
there is a neighborhood $V$ of $p$, relatively compact in $V_1$ 
and a positive constant $s<1$,
independent of $q$, such that for every $f \in \mathcal
O_{J}(\Delta,D)$ satisfying $f(0) \in D \cap V$ we have $f(s\Delta)
\subset D \cap U$. This gives the inequality (\ref{e2}). \qed

\section{Scaling and  estimates of the Kobayashi-Royden metric} 
In this Section we present a precise lower estimate on the Kobayashi-Royden
infinitesimal metric on a strictly pseudoconvex domain in
$(M,J)$.

\begin{theorem}\label{THM}
Let $M$ be a real $2n$-dimensional  manifold with an almost complex
structure $J$ and
let $D=\{\rho<0\}$ be a relatively compact domain in $(M,J)$.
We assume that $\rho$ is a $\CC^2$ defining function of $D$,
strictly $J$-plurisubharmonic in a neighborhood of $\bar{D}$. Then
there exists a positive constant $c$  such that~:

\begin{equation}\label{e3}
K_{(D,J)}(p,v) \geq c\left[\frac{|\partial_J\rho(p)(v - iJ(p)v)|^2}
{|\rho(p)|^2} + 
\frac{\|v\|^2}{|\rho(p)|}\right]^{1/2},
\end{equation}
for every $p \in D$ and every $v \in T_pM$.
\end{theorem}

We start with the small almost complex deformations of the standard
structure. In the second subsection, we consider the case of an
arbitrary almost complex structure, not necessarily close to the
standard one. We use non-isotropic dilations in special coordinates
``reducing'' an almost complex structure in order to represent a
strictly pseudoconvex hypersurface on an almost complex manifold 
as the Siegel sphere equipped with an arbitrary
small deformation of the standard structure. We stress
that such a representation cannot be obtained by the isotropic
dilations of Lemma 1 since the limit hypersurface is just a
hyperplane. 

\subsection{Small deformations of the standard structure}

We start the proof of Theorem~\ref{THM} with the following~:

\begin{proposition}\label{thm2}
Let $D=\{\rho < 0\}$ be a bounded domain in $\C^n$, where $\rho$ is a
 $\CC^2$ defining function of $D$, strictly $J_{st}$-plurisubharmonic in
 a neighborhood of $\bar{D}$. Then there exist positive constants $c$
 and $\lambda_0$ such that for every almost complex structure $J$
 defined in a neighborhood of $\bar{D}$ and such that
 $\|J-J_{st}\|_{\CC^2(\bar{D})} \leq \lambda_0$ estimate~(\ref{e3}) is
satisfied for every $p \in D,$ $v \in \C^n$.
\end{proposition}

\vskip 0,1cm
\noindent{\it Proof}. We note that according to Proposition~\ref{thm3}
 (see estimate (\ref{localhyp}))
it is sufficient to prove the inequality 
near $\partial D$. Suppose by contradiction that there exists a
sequence $(p^{\nu})$ of points in $D$ converging to a boundary point
$q$, a sequence $(v^{\nu})$ of unitary vectors and a sequence $(J_\nu)$ 
of almost complex structures defined in a neighborhood of $\bar{D}$, 
satisfying 
$\|J_\nu-J_{st}\|_{\CC^2(\bar{D})}\rightarrow_{\nu \rightarrow \infty} 0$, 
such that the quotient

\begin{equation}
\label{quot1}
K_{(D\cap U,J_\nu)}(p^{\nu},v^{\nu})\left[\frac{|\partial_{J_{\nu}} 
\rho(p^{\nu})(v^{\nu} - iJ_{\nu}(p^{\nu})v^{\nu})|^2}{|\rho(p^{\nu})|^2} 
+ \frac{\|v^{\nu}\|^2}{|\rho(p^{\nu})|}\right]^{-1/2}
\end{equation}
tends to $0$ as $\nu$ tends to $\infty$, where
$U$ is a neighborhood of $q$.
For sufficiently large $\nu$ denote by
$\delta_{\nu}$ the euclidean distance from $p^{\nu}$ to the
boundary of $D$ and by $q^{\nu} \in \partial D$ the unique point such
that $\vert p^{\nu} - q^{\nu} \vert = \delta_{\nu}$. Without loss of
generality we assume that $q = 0$, that $T_0(\partial D) = \{z:=('z,z_n) \in
\C^n : Re(z_n) = 0\}$ and that $J_\nu(q^\nu) = J_{st}$ for every $\nu$.

Consider a sequence of biholomorphic (for the standard structure)
transformations $T^{\nu}$ in a neigborhood of the origin, such that
$T^\nu(q^{\nu}) = 0$ and such that the image $D^{\nu} : =T^\nu(D)$
satisfies
$$
T_0(\partial D^{\nu})=\{ z \in \C^n : Re(z_n) = 0\}.
$$
We point out that the sequence $(T^{\nu})_\nu$ converges uniformly to the
identity map since $q^{\nu} \rightarrow q=0$ as $\nu \rightarrow \infty$
and hence that the sequence $((T^\nu)^{-1})_\nu$ is bounded. We
still denote by $J_\nu$ the direct image
$(T^\nu)_*(J_\nu)$. Let $U_1$ be a neighborhood of the origin such that
$\bar{U} \subset U_1$. For sufficiently large $\nu$ we have 
$T^\nu(U) \subset U_1$. We may assume that every domain $D^{\nu}$ is
defined on $U_1$ by
$$
D^\nu \cap U_1 = \{z \in U_1 : \rho^{\nu}(z) := 
Re(z_n) + |'z|^2 +\mathcal O(|z|^3) <0\},
$$ 
and that the
sequence $(\hat p^\nu = T^{\nu}(p^{\nu}) =(0',-\delta_\nu))_\nu$ is on
the real inward normal to $\partial D^{\nu}$ at 0. Of course, the
functions $\rho^{\nu}$ converge uniformly with all derivatives to the
defining function $\rho$ of $D$. In what follows we omit the hat and
write $p^{\nu}$ instead of $\hat{p}^{\nu}$.

Denote by $R$ the function
$$
{R}(z)=Re(z_n) + |'z|^2 + (Re(z_n) + \vert 'z \vert^2)^2.
$$
There is a neighborhood $V_0$ of the origin in $\mathbb C^n$
such that the function
${R}$ is strictly $J_{st}$-plurisubharmonic on $V_0$. Fix $\alpha > 0$
small enough 
such that the point $z^\alpha=('0,-\alpha)$ belongs to $V_0$.
Consider the dilation $\Lambda_\nu$ defined on
$\C^n$ by $\Lambda_\nu(z) =
({(\alpha / \delta_\nu)^{1/2}}'z,(\alpha/\delta_\nu)z_n)$. 
If we set $J^\nu :=\Lambda_\nu \circ J_\nu \circ (\Lambda_\nu)^{-1}$ 
then we have~: 

\begin{lemma}\label{Lemma1}
 $\lim_{\nu \rightarrow \infty}J^\nu = J_{st}$, uniformly on compact 
subsets of $\C^n$. 
\end{lemma}

\noindent{\it Proof}. Considering $J$ as a matrix valued function,
we may assume that the Taylor expansion of $J_\nu$ at the origin is given by 
$J_\nu = J_{st} + L_\nu(z) + \mathcal O(|z|^2)$ on $U$, uniformly with 
respect to $\nu$. 
Hence ${J}^\nu(z^0)(v) = J_{st}(v) +
{L}_\nu('z,(\delta_\nu/\alpha)^{1/2}z_n)(v)
+ \mathcal O(|(\delta_\nu|) 
\ \|v\|$.
Since $\lim_{\nu \rightarrow \infty}L_\nu = 0$ by assumption, 
we obtain the desired result. \qed

\vskip 0,1cm
Let $\tilde{\rho}^\nu:=(\alpha / \delta_\nu) \rho^\nu \circ 
\Lambda_\nu^{-1}$ and
$G^\nu:=\{z \in \Lambda_\nu(U_1) : \tilde{\rho}^\nu(z) < 0\}$.
Then the function $R^\nu:= \tilde{\rho}^\nu + (\tilde{\rho}^\nu)^2$
converges with all its derivatives to ${R}$, 
uniformly on compact subsets of $\mathbb C^n$.
Hence $R^\nu$ is strictly plurisubharmonic on $V_0$
and according to Lemma~\ref{Lemma1} there is a positive constant $C$ 
such that for sufficiently large $\nu$ the function $R^\nu - C|z|^2$ 
is strictly $J^\nu$-plurisubharmonic on $V_0$.
Since $\sup_{z \in G^\nu \cap \partial V_0} (R^\nu(z) -C|z|^2) =-C'<0$, 
the function 
$$
\tilde{R}^\nu:=\left\{
\begin{array}{lll}
R^\nu - C|z|^2 & {\rm on} & D^\nu \cap V_0\\
 & & \\
-\frac{C'}{2} & {\rm on} & D^\nu \backslash V_0
\end{array}
\right.
$$
is $J^\nu$-plurisubharmonic on $G^\nu$, strictly $J^\nu$-plurisubharmonic
on $G^\nu \cap V_0$. 
Since $z^\alpha$ belongs to $V_0$, it follows from the Proposition~\ref{thm3} 
(see estimate (\ref{localhyp})) that there exists a  positive constant
$C' > 0$ such that for sufficiently large $\nu$ we have~:

$$
K_{(G^\nu,J^\nu)}(z^\alpha,v) \geq  C'\|v\|
$$
for every $v \in \C^n$. 

Moreover for $v \in \C^n$ and for sufficiently large $\nu$ we have~:
\begin{eqnarray*}
& &   
K_{(D^\nu\cap U_1, J_\nu)}(p^\nu,v) =
K_{(G^\nu,J^\nu)}(z^\alpha,\Lambda_\nu(v)) \geq  C'\parallel \Lambda_\nu(v)
\parallel.
\end{eqnarray*}
This gives the inequality~:
$$
K_{(D^{\nu} \cap U,J_\nu)}(p^\nu,v) \geq C' \left(
\frac{\alpha |v_1|^2}{\delta_\nu} + \cdots + 
\frac{\alpha |v_{n-1}|^2}{\delta_\nu} + 
\frac{\alpha^2 |v_n|^2}{\delta_\nu^2}\right)^{1/2}.
$$
Since $C_1\delta_\nu$ is equivalent to $|\rho(p^\nu)|$ as $\nu
\rightarrow \infty$, we obtain that there
is a positive constant $C''$ such that
$$
K_{(D^{\nu} \cap U,J_\nu)}( p^\nu,v) \geq C'' \left(
\frac{\|v\|^2}{|\rho(p^\nu)|} + 
\frac{|v_n|^2}{|\rho(p^\nu)|^2}\right)^{1/2}.
$$
Since $J_{\nu}(0) = J_{st}$, we have  $|\partial\rho(p^\nu)(v - 
iJ_\nu(p^\nu)(v))|^2 = |\partial_{J_{st}}\rho(p^\nu)(v)|^2 +
\mathcal O(\delta_{\nu})\parallel v \parallel^2 = \vert v_n \vert^2 
 + \mathcal O(\delta_\nu)\parallel v \parallel^2$.
Hence there exists a positive constant $\tilde{C}$ such that 
$$
K_{(D^{\nu} \cap U,J_\nu)}( p^\nu,v) \geq \tilde{C} \left(
\frac{\|v\|^2}{|\rho(p^\nu)|} + 
\frac{|\partial_J\rho(p^\nu)(v-iJ_\nu(p^\nu)(v))|^2}
{|\rho(p^\nu)|^2}\right)^{1/2},
$$
contradicting the assumption on the quotient 
(\ref{quot1}). This proves the desired estimate. \qed

We  have the following corollary~:
\begin{corollary}\label{cor3}
Let $(M,J)$ be an almost complex manifold. Then every $p \in M$ has a 
basis of complete hyperbolic neighborhoods.
\end{corollary}

\proof Let $p \in M$. According to Example~1 there exist a
neighborhood $U$ of $p$ and a diffeomorphism $z:U \rightarrow \B$,
centered at $p$, such that the function $|z|^2$ is strictly
$J$-plurisubharmonic on $U$ and 
$\|z_\star(J)-J_{st}\|_{\mathcal C^2(U)} \leq \lambda_0$. 
Hence the open ball $\{x \in \mathbb C^n :
\|x\|<1/2\}$ equipped with the structure  $z_\star(J)$ satisfies the
hypothesis of Theorem \ref{thm2}. Now the estimate on the
Kobaysahi-Royden metric given by this theorem implies that this ball
is  complete hyperbolic by the standard integration argument. 
\qed  

\subsection{Arbitrary  almost complex structures}

We turn now to the proof of Theorem~\ref{THM} on an
arbitrary  strictly pseudoconvex domain in an almost
complex manifold $(M,J)$ ($J$ is not supposed to
be a small deformation of the standard structure).
In view of Proposition~\ref{thm3} it suffices to prove
the statement in a neighborhood $U$ of  a boundary point  $q \in \partial D$.
Considering local coordinates $z$ centered at $q$, we may assume that
$D \cap U$ is a domain in $\mathbb C^n$ and
$0 \in \partial D$, $J(0) = J_{st}$.
The idea of the proof is to reduce the situation to the case of a small
deformation of the standard structure considered in Proposition~\ref{thm2}.
In the case of real dimension four Theorem~\ref{THM} is a direct corollary of
Proposition~\ref{thm2}. In the case of arbitrary dimension the proof of
Theorem~\ref{THM} requires a slight modification of Proposition~\ref{thm2}.
So we treat this case seperately.

\subsubsection{Case where $dim M = 4$}
According to \cite{sik94} Corollary~3.1.2,
there exist a neighborhood $U$ of $q$ in $M$ and complex coordinates
$z=(z_1,z_2) : U \rightarrow \mathbb B_2 \subset \mathbb C^2$, $z(0) =
0$ such that $z_*(J)(0) = J_{st}$ and moreover, a map $f: \Delta
\rightarrow \mathbb B$ is $J':= z_*(J)$-holomorphic if it satisfies the
equations 

\begin{eqnarray}
\label{Jhol}
\frac{\partial f_j}{\partial \bar \zeta} =
A_j(f_1,f_2)\overline{\left ( \frac{\partial f_j}
{\partial \zeta}\right ) }, j=1,2
\end{eqnarray} 
where $A_j(z) =  O(\vert
z \vert)$, $j=1,2$.

In order to obtain such coordinates, one can consider two transversal
foliations of $\mathbb B$ by $J'$-holomorphic curves (see \cite{nw63})
and then take these curves into the lines $z_j = const$ by a local
diffeomorphism. The direct image of the almost complex structure $J$
under such a diffeomorphism has a diagonal matrix $ J'(z_1,z_2) =
(a_{jk}(z))_{jk}$ with $a_{12}=a_{21}=0$ and $a_{jj}=i+\alpha_{jj}$
where $\alpha_{jj}(z)=\mathcal O(|z|)$ for $j=1,2$.
We point out that the lines $z_j = const$ are
$J$-holomorphic after a suitable parametrization (which, in general,
is not linear). 

In what follows we omit the prime and denote this structure again by
$J$. We may assume that the complex tangent space $T_0(\partial D)
\cap J(0) T_0(\partial D) = T_0(\partial D) \cap i T_0(\partial D)$ is
given by $\{ z_2 = 0 \}$.
In particular, we have the following expansion for the defining
function $\rho$ of $D$ on $U$~:
$\rho(z,\bar{z}) = 2 Re(z_2) + 2Re K(z) + H(z) + \mathcal O(\vert z
\vert^3)$, where
$K(z)  = \sum k_{\nu\mu} z_{\nu}{z}_{\mu}$, $k_{\nu\mu} =
k_{\mu\nu}$ and 
$H(z) = \sum h_{\nu\mu} z_{\nu}\bar z_{\mu}$, $h_{\nu\mu} =
\bar h_{\mu\nu}$.

\begin{lemma}\label{PP}
The domain $D$  is strictly $J_{st}$-pseudoconvex near the origin.
\end{lemma}

\noindent{\sl Proof of Lemma~\ref{PP}}. Consider a complex vector
 $v=(v_1,0)$ tangent to $\partial D$ at the origin.  Let $f:\Delta
 \rightarrow \mathbb C^2$ be a $J$-holomorphic disc centered at the
 origin and tangent to $v$: $f(\zeta) = v\zeta + \mathcal O(\vert
 \zeta \vert^2)$.  Since $A_2 = \mathcal O(\vert z \vert)$, it follows
 from the $J$-holomorphy equation (\ref{Jhol}) that
 $(f_2)_{\zeta\bar\zeta}(0) = 0$. This implies that $(\rho \circ
 f)_{\zeta\bar\zeta}(0) = H(v).$ Thus, the Levi form with respect to
 $J$ coincides with the Levi form with respect to $J_{st}$ on the
 complex tangent space of $\partial D$ at the origin.\qed

\vskip 0,1cm
Consider the non-isotropic dilations $\Lambda_{\delta}: (z_1,z_2) \mapsto
(\delta^{-1/2}z_1,\delta^{-1}z_2) = (w_1,w_2)$ with $\delta > 0$. 
If $J$ has the above
diagonal form in the coordinates $(z_1,z_2)$ in $\mathbb C^2$, then
its direct image  $J_{\delta}= (\Lambda_{\delta})_*(J)$ has the form
$J_{\delta}(w_1,w_2) =(a_{jk}(\delta^{1/2}w_1,\delta w_2))_{jk}$
and so $J_{\delta}$ tends to $J_{st}$ in the $\mathcal C^2$ norm as $\delta
\rightarrow 0$. On the other hand, $\partial D$ is, in the coordinates
$w$, the zero set of the function 
$\rho_{\delta}= \delta^{-1}(\rho \circ \Lambda_{\delta}^{-1})$.
As $\delta \rightarrow 0$, the function $\rho_{\delta}$ tends to 
the function $2 Re w_2 + 2 Re K(w_1,0) + H(w_1,0)$ which defines a
$J_{st}$- strictly pseudoconvex domain by Lemma~\ref{PP}.
So we may apply Proposition~\ref{thm2}.
This proves Theorem~\ref{THM} in dimension 4.

\subsubsection{Case where $dim M = 2n$.} 
In this case  Proposition \ref{thm2} is not directly applicable since
$J$ can not be deformed by the non-isotropic dilations to the standard
structure. Instead we use the invariance of the Levi form with respect
to the non-isotropic dilations.

We suppose that in a neighborhhod of the origin we have $J = J_{st} +
{\mathcal O}(\vert z \vert)$.
We also may assume that in these coordinates the defining function
$\rho$ of $D$ has the form $\rho = 2Re z_n + 2ReK(z) + H(z) + 
\mathcal O(\vert z \vert^3)$, where $K$ and $H$ are defined similarly
to the 4-dimensional case and $\rho$ is strictly $J$-plurisubharmonic
at the origin. We use the notation $z = ('z,z_n)$.

Consider the non-isotropic dilations $\Lambda_{\delta} : ('z,z_n)
\mapsto (w',w_n) = {(\delta^{-1/2}}'z,\delta^{-1}z_n)$ and set
$J_{\delta} = (\Lambda_{\delta})_*(J)$. Then $J_{\delta}$ tends to the
almost complex structure $J_0(z)= J_{st} + L('z,0)$ where 
$L('z,0) = (L_{kj}('z,0))_{kj}$
denotes a matrix with $L_{kj} = 0$ for $k = 1,...,n-1$, $j = 1,...,n$,
$L_{nn} = 0$ and $L_{nj}('z,0)$, $j=1,...,n-1$ being (real) linear
forms in $'z$.  

Let $\rho_{\delta} = \delta^{-1}(\rho \circ \Lambda_{\delta}^{-1})$.
As $\delta \rightarrow 0$, the function $\rho_{\delta}$ tends to 
the function $\tilde{\rho} = 2 Re z_n + 2 Re K('z,0) + H('z,0)$ in the
$\mathcal C^2$ norm. By the invariance of the Levi form we have
${\mathcal L}^J(\rho)(0)(\Lambda_{\delta}^{-1}(v)) = {\mathcal
  L}^{J_\delta}(\rho \circ \Lambda_{\delta}^{-1})(0)(v)$. Since $\rho$ is
strictly $J$-plurisubharmonic, multiplying by $\delta^{-1}$ and
passing to the limit at the right side  as $\delta \longrightarrow 0$ , 
we obtain that
${\mathcal L}^{J_0}(\tilde \rho)(0)(v) \geq 0$ for any $v$. Now let $v =
('v,0)$. Then $\Lambda_{\delta}^{-1}(v) = \delta^{1/2}v$ and so 
${\mathcal L}^J(\rho)(0)(v) = {\mathcal
  L}^{J_\delta}(\rho_{\delta})(0)(v)$. Passing to the limit as $\delta$
tends to zero, we obtain that  ${\mathcal L}^{J_0}(\tilde \rho)(0)(v) > 0$
for any $v = ('v,0)$ with $'v \neq 0$.

Consider now the function $R=\tilde{\rho} + \tilde{\rho}^2$. 
Then ${\mathcal L}^{J_0}(R)(0)(v) = {\mathcal L}^{J_0}(\tilde \rho)(0)(v)
+ 2 v_n \overline v_n$, so $R$ is strictly $J_0$-plurisubharmonic in a
neighborhood of the origin.
Thus the functions $R^{\nu}$ used in the proof of
Proposition~\ref{thm2} are strictly $J^\nu$-plurisubharmonic and their Levi
forms are bounded from below by a positive constant independent of $\nu$.
This allows to use Proposition~\ref{thm3} and the proof can be
proceeded quite similarly to the proof of Proposition~\ref{thm2}
without any changes.  \qed

\begin{remark} As it was brought to our attention by J.P.Rosay, our proof
of Theorem~\ref{THM} also gives the following local version. If $D$ is an
arbitrary domain (not necessarily relatively compact) in an almost complex
manifold $(M,J)$, strictly pseudoconvex at a boundary point $q$, then the
estimate~(\ref{e3}) still holds for any point $p\in D$ in a neighborhood
of $q$.
\end{remark}

\end{document}